\documentclass[%
 reprint,
 amsmath,amssymb,
 aps,
]{revtex4-2}

\usepackage{graphicx}
\usepackage{dcolumn}
\usepackage{bm,gensymb,multirow}

\begin{document}

\preprint{APS/123-QED}

\title{Detecting self-organising patterns in crowd motion: Effect of optimisation algorithms}

\author{Samson Worku}
\author{Pratik Mullick}%
 \email{pratik.mullick@pwr.edu.pl}
\affiliation{%
 Department of Operations Research and Business Intelligence, Wrocław University of Science and Technology, Wrocław, Poland
}%

\begin{abstract}
The escalating process of urbanization has raised concerns about incidents arising from overcrowding, necessitating a deep understanding of large human crowd behavior and the development of effective crowd management strategies. This study employs computational methods to analyze real-world crowd behaviors, emphasizing self-organizing patterns. Notably, the intersection of two streams of individuals triggers the spontaneous emergence of striped patterns, validated through both simulations and live human experiments. Addressing a gap in computational methods for studying these patterns, previous research utilized the pattern-matching technique, employing the Nelder-Mead Simplex algorithm for fitting a two-dimensional sinusoidal function to pedestrian coordinates. This paper advances the pattern-matching procedure by introducing Simulated Annealing as the optimization algorithm and employing a two-dimensional square wave for data fitting. The amalgamation of Simulated Annealing and the square wave significantly enhances pattern fitting quality, validated through statistical hypothesis tests. The study concludes by outlining potential applications of this method across diverse scenarios.
\end{abstract}

\maketitle


\onecolumngrid

\section{Introduction}
\label{intro}

\begin{figure}[h!]
    \centering
    \includegraphics[width=\textwidth]{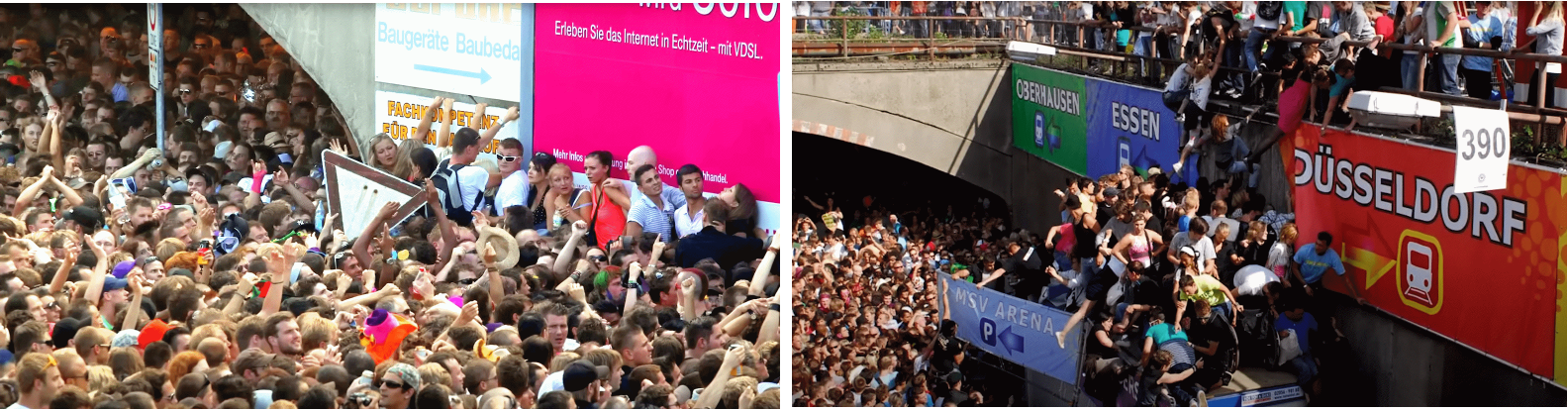}
    \caption{Crowd disaster in Duisburg, Germany during Love Parade music festival in 2010. The crowd trapped under the tunnel tried to escape, which later resulted in stampedes. Images taken as screenshots at 11:10 (left) and 12:54 (right) from https://www.youtube.com/watch?v=rq5bcWOjdrU}.
    \label{love_parade}
\end{figure}

Gathering of a large number of people is generally known as a crowd. With the increased rate of urbanisation, the world has experienced several catastrophic incidents with death tolls, resulting from overcrowding at social gatherings. The death toll was 21 in Love Parade music festival in Duisburg, Germany in 2010 (Fig. \ref{love_parade}), 3 during a stampede in UEFA Champions League Final in Turin, Italy (with 1500+ injured) in 2021, 150 in Halloween celebrations in Seol, South Korea in 2022 and unfortunately the list goes on. The deadliest crowd accidents belong to the annual Hajj pilgrimage in Mecca and Mina, Saudi Arabia, killing 5000+ people in a span of about 25 years. Most of these incidents are consequences of poor crowd management by the organisers.

Research on the motion of a crowd finds its importance in providing efficient techniques for crowd management. In urban spaces, in most of the cases, the crowd is basically a group of pedestrians walking through the streets, traffic-signal junctions, railway or metro stations, shopping malls, sports stadiums etc. During emergency situations (such as fire, earthquakes, sudden appearance of a gunman etc.), successful evacuation of the crowd needs well-researched preventive measures. However, even without an emergency escape situation, managing a large crowd is definitely a challenging job, e.g. when a large number of people come out of a stadium after a sports event.

When exposed to specific settings, individuals in a crowd have a `natural' method of organising and structuring themselves. This self-organisation in crowds enables researchers to investigate and comprehend its underlying behaviour, that is consistent with the crowd's `natural' behaviour. Self-organisation, in the context of crowd management, refers to the spontaneous formation of organised patterns by pedestrians without any external guidance or control. It is the natural tendency of individuals within a crowd to align themselves in a unidirectional manner, resulting in the formation of `lanes'. The formation of spatio-temporal patterns in a system as a result of socio-physical interactions of human beings occur without any external design or structure.

\begin{figure}[h!]
    \centering
    \includegraphics[width=\textwidth]{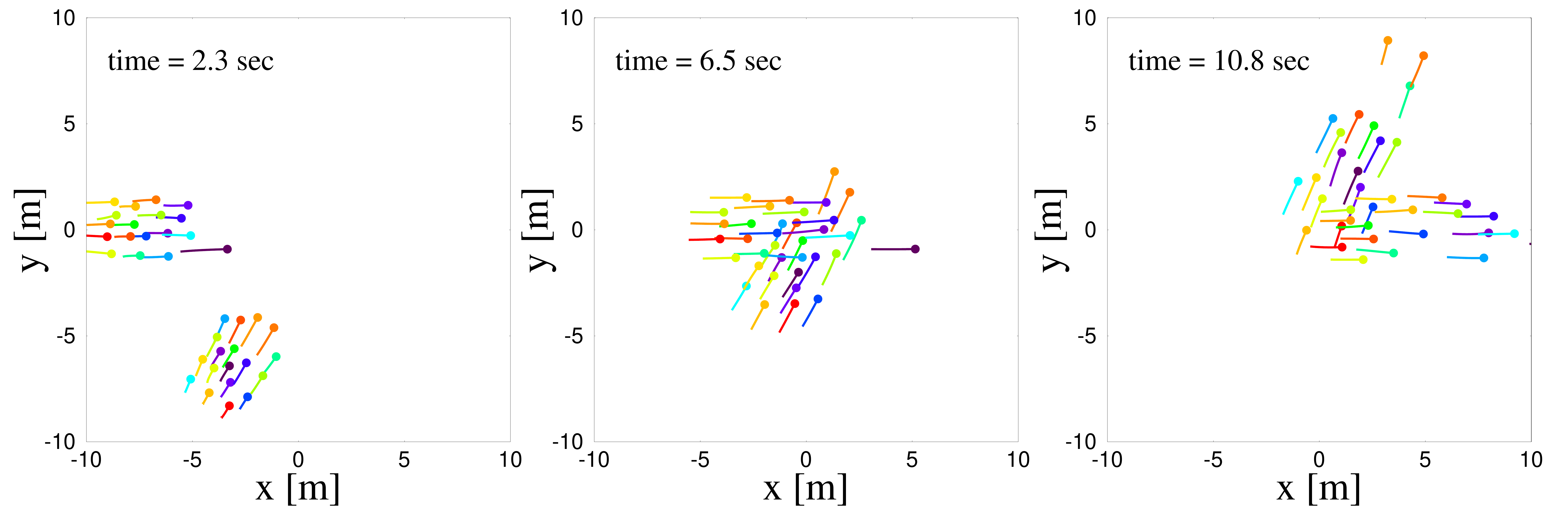}
    \caption{Typical experimental data for crossing flows of pedestrians on a 2D surface ($x-y$ plane) at a crossing angle = 60\degree. The pedestrians are denoted by dots and the tails behind them are the distances travelled by the pedestrians in previous 1.25 sec. The data is published in \cite{mullick_ploscb}.}
    \label{traces_cross_flow}
\end{figure}

To create effective measures of crowd management, one needs to understand the behaviour of a crowd \cite{haghani2018} in a wide variety of situations. Especially in an urban set-up, where there could be several types of emergency-escape situations \cite{helbing2000}, crowd management strategies become more essential for crowd evacuation \cite{thompson1995,xie2021}. The critical factor is the crowd density rather than the crowd size. At high densities, the behavior of a human crowd has several analogies to that of a flowing fluid. The most common crowd scenarios that could be seen in urban areas are uni-directional flows \cite{seyfried2005}, bidirectional flows \cite{sharifi2020}, multi-directional flows \cite{cao2017}, crossing flows \cite{mullick_ploscb,zanlungo1,zanlungo2}, flow through a T-shaped junction \cite{zhang2011}, single pedestrian crossing through a dense static crowd \cite{nicolas2019}, bottlenecks \cite{hoog2005,seyfried2009,nicolas2017} etc. Each of these crowd situations have unique self-organising patterns. For example, formation of lanes are seen when two groups of people try to cross through each other in a sidewalk \cite{social_force}. This particular crowd situation is known as a counterflow, which is a special case of crossing flows when the crossing angle is $180\degree$. In Fig. \ref{traces_cross_flow} a typical human data for pedestrian flows crossing at an angle of $60\degree$ are shown.

To understand the crowd behaviour one has to rely on experimental data using real human beings, but conducting these experiments is very difficult from financial and ethical points of view. For this reason we turn to model the crowd using available experimental data, thanks to the improvement of computational facilities in the past few decades. The research area of crowd simulation concerns the design and use of simulation algorithms to understand, predict and reproduce the behaviour of human-crowds \cite{wouter_julien,saeed2022}. The immediate output of a crowd simulation algorithm are the time-dependent trajectories of humans, as realistic as possible. Studies of collective motion have also been performed for several other social groups of the animal kingdom \cite{john2004,gau2012,giu2015,totzeck2020}.

\begin{figure}[h!]
    \centering
    \includegraphics[width=\textwidth]{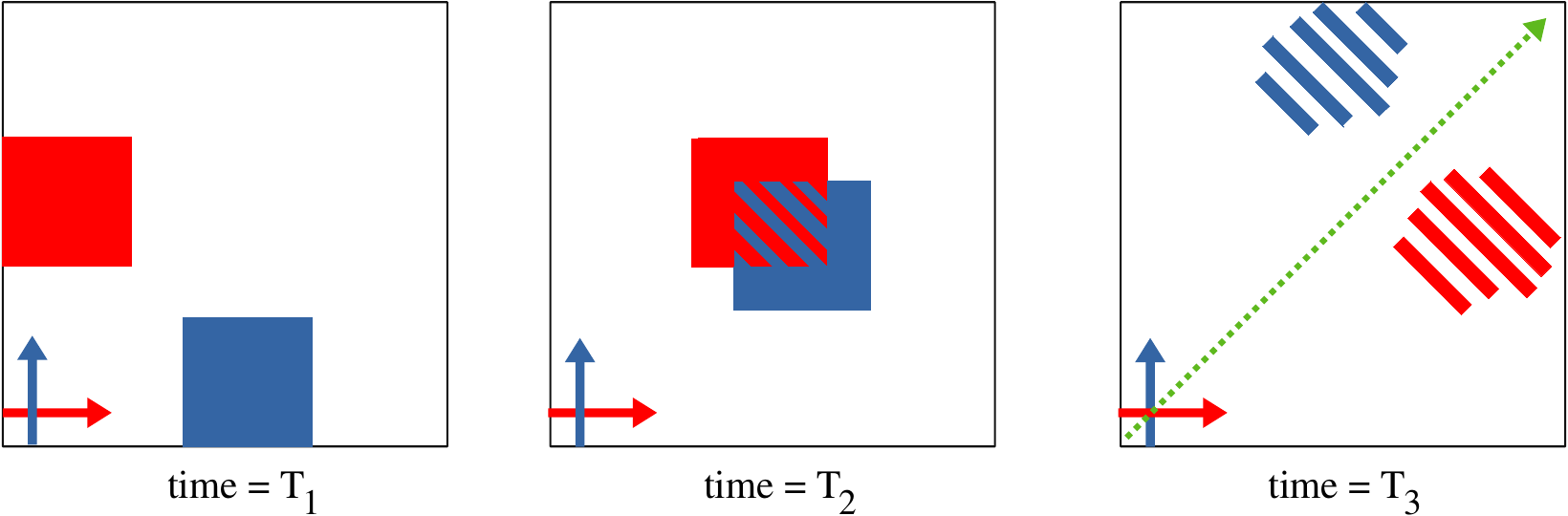}
    \caption{Schematic representation for the formation of stripes when two groups of people cross through each other. The figure has been shown for three instances - before crossing ($T_1$), during crossing ($T_2$) and after crossing ($T_3$); and hence $T_1<T_2<T_3$. The two groups before crossing are denoted by blue and red squares, whose direction of motion is denoted by arrows of the same color. The green dotted arrow denotes the bisector of the crossing angle.}
    \label{stripe_schematic}
\end{figure}

Studies of human-crowd dynamics by both experiments and real-time simulations have diverse range of applications from entertainment, such as computer games and movies, to safety-critical analysis, i.e., to improve pedestrian traffic flow and prevention of crowd disasters \cite{thompson1995,helbing2000,springer2002,bohannon2005,enguica2023}. Pedestrian traffic flow has been studied empirically in a wide variety of situations, using both experimental methods and motion tracking of real crowds \cite{helbing2007,johansson2008}. When two streams of pedestrians cross at an angle, striped patterns spontaneously emerge as a result of local pedestrian interactions. This phenomenon is schematically illustrated in Figure \ref{stripe_schematic}. Several urban situations produce crossing flows, such as streams of pedestrians crossing at a sidewalk intersection, or subway commuters passing each other when entering and exiting a metro car. It is very common to notice that pedestrians in a crosswalk often form multiple lanes of traffic. Such spontaneous pattern formation is an example of self-organised collective behaviour, a topic of intense interdisciplinary interest \cite{helbing2005,mehdi2012}.

In the early stages of human-crowd research back in 1977, an empirical hypothesis posited that striped patterns in crossing flows should align perpendicularly to the bisector of the crossing angle \cite{naka1977}. However, since that early proposition, there has been a scarcity of experimental evidence supporting this phenomenon. While simulation-based studies have explored crossing flows of pedestrians, none have definitively validated the `bisector-normal hypothesis.' A stripe represents a traveling wave that aligns with the mean direction of the two flows, facilitating forward movement for pedestrians within it \cite{helbing2005}. This striped pattern is crucial for optimizing pedestrian flow, minimizing collision-avoidance situations, and consequently enhancing the average walking speed. Although a handful of subsequent human studies have delved into oblique crossing angles \cite{guo,wong,aghabayk}, none have specifically analyzed the presence or properties of stripe patterns. Notably, in \cite{mullick_ploscb}, pioneering computational methods were introduced to detect the existence of striped patterns, paving the way for a comprehensive study of their geometric properties. The outcomes from these methods confirmed that the observed striped patterns align perpendicularly to the bisector of the crossing angle. This study represents a significant stride in substantiating and exploring the long-standing hypothesis regarding the directional alignment of striped patterns in crossing flows.

In this paper we consider the experimental data of crossing flows of pedestrians \cite{mullick_ploscb}, where two groups of volunteers crossed each other at 7 different crossing angles from $0\degree$ to $180\degree$, at intervals of $30\degree$. Emergence of self-organising `striped' patterns has already been demonstrated even for small groups pedestrians crossing each other. In one of the methods developed in \cite{mullick_ploscb}, viz. the pattern-matching technique, the expected pattern of alternate and parallel stripes in the crossing region (see Fig. \ref{stripe_schematic}) was captured by employing a two-dimensional sinusoidal function for fitting pedestrian coordinates. The fitting was done using the Nelder-Mead simplex algorithm which maximises an objective function (see Sect.~\ref{method}) designed accurately to capture the striped pattern.

The primary objective of this paper is to provide improvements in the pattern fitting procedure. Our essential approach is two-fold: (i) to use Simulated Annealing as the optimisation algorithm and (ii) to employ a two-dimensional square wave to capture the striped pattern. Our findings indicate that the changes in the optimisation strategy in the two above ways indeed results in an improved fitting procedure, which was verified using statistical hypothesis tests. Insights gained from this research could be useful for creating effective crowd simulation models. The remaining part of this paper is organised as follows: in the next section, viz. Sect. \ref{method} we describe the methodological tools and techniques used in this research, followed by results and discussion in Sect. \ref{res.dis}. Finally in Sect. \ref{con} we make some concluding remarks.

\section{Materials and Methods}
\label{method}

The results presented in this paper is based on experimental data that was collected by employing volunteer pedestrians and recording their trajectories. Using these data we built computational techniques to detect the presence of striped patterns. In this section we briefly describe the experimental details, followed by the computational strategies utilised for this research.

\subsection{Experimental details}
\label{expt}



The experimental data for the pedestrian crossing flows utilized in this study were acquired through experimental trials conducted on the University of Rennes campus in France, as detailed in prior works \cite{pedinteract_cecile,mullick_ploscb}. This dataset is publicly accessible via the repository at \url{https://doi.org/10.5281/zenodo.5718430}. Two distinct groups of volunteer participants, comprising 36 individuals on Day 1 and 38 on Day 2, were randomly divided into two groups, with each group containing 18 or 19 individuals. These participants were given instructions to traverse a sports hall to reach the opposite side. Initial positions were prearranged to necessitate crossings at seven different angles (ranging from $0\degree$ to $180\degree$ in $30\degree$ increments). Throughout each trial, we captured the head trajectories of all pedestrians as time series data at a frequency of 120 Hz using an infrared camera-based motion capture system (VICON). A total of 116 trials were recorded across all crossing angles, with 10 trials corresponding to a $0\degree$ crossing angle. Consequently, we retained 106 trials for our research, focusing on the presence of striped patterns in the crossing region, a phenomenon that arises when the two groups walk in different directions and intersect.

Subsequently, the acquired data underwent a low-pass filtering process to attenuate oscillations resulting from the natural gait movements of walking pedestrians. Specifically, we applied a forward-backward 4th-order Butterworth filter with a cut-off frequency set at 0.5 Hz. Figure \ref{traces_cross_flow} illustrates the trajectories of all pedestrians in a representative trial using the filtered trajectory data.

\subsection{Pattern-matching technique}
\label{pat.match}

When two streams of people cross at oblique incidence we expect striped patterns as demonstrated in Figure \ref{stripe_schematic}. However, when we merely plot the experimental data, which is the trajectories of all the participants, the expected patterns are not obvious to the eyes, as in Figure \ref{traces_cross_flow}. For this purpose a pattern-matching technique was developed \cite{mullick_ploscb}, which we describe here.

To capture the expected pattern of alternate and parallel stripes in the crossing region, a two-dimensional sinusoidal function was employed to fit pedestrian coordinates $(x,y)$. The sinusoid $f$ is defined as:
\begin{equation}
    f(x,y;\gamma,\lambda,\psi) = \sin\left(\frac{2\pi X}{\lambda}+\psi\right),
    \label{gabor_sine}
\end{equation} where:
\begin{equation}
X = x\sin\gamma-y\cos\gamma,
\end{equation} rotates the coordinates using angle $\gamma$ to represent the orientation of stripes with respect to the $x$-axis. The parameter $\lambda$ represents the wavelength of the sine curve, corresponding to the spatial separation between stripes from the same group, and $\psi$ denotes a phase offset. The fitting procedure involved maximising the function $C$:
\begin{equation}
C = \sum_{\text{group 1}}\frac{f(x,y)}{N_1} + \sum_{\text{group 2}}\frac{-f(x,y)}{N_2},
\label{c.eq}
\end{equation} where $N_1$ and $N_2$ are the number of pedestrians in the two groups, respectively. The maximisation of $C$ was performed using the Nelder-Mead simplex algorithm. The maximum possible value of $C$ is 2, indicating an ideal fit when pedestrian positions from both groups precisely match the crests and troughs of the sinusoid. By maximising $C$ through fitting $f$ to pedestrian positions, we obtain the stripe orientation $\gamma$ and spatial separation $\lambda$.

\begin{figure}[h!]
    \centering
    \includegraphics[width=\textwidth]{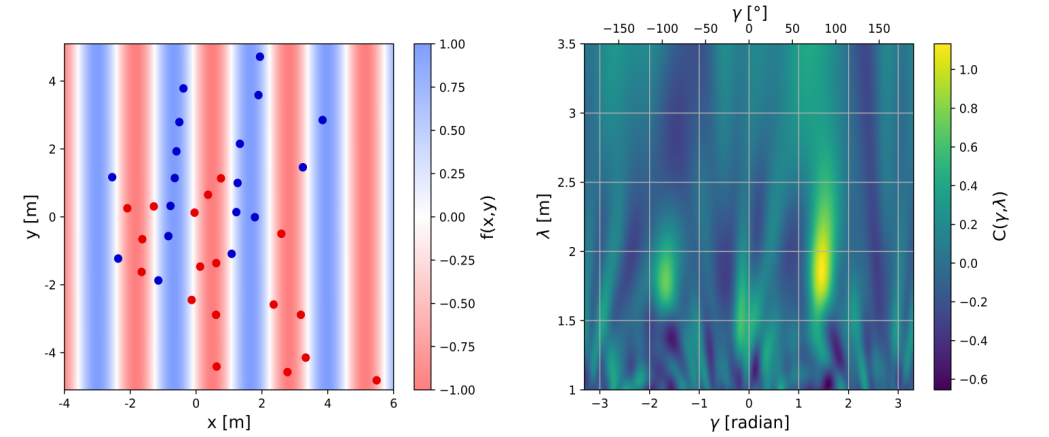}
    \caption{Demonstration of pattern-matching technique for a typical trial. The left panel shows the 2D sine wave fitted on pedestrian coordinates from the two groups, denoted by blue and red dots, and the right panel shows the corresponding objective function $C$ as a function of $\gamma$ and $\lambda$.}
    \label{pat.match.sine}
\end{figure}

It is worth mentioning here that the idea of fitting a two dimensional sinusoid (Equation \ref{gabor_sine}) to capture the periodic and parallel striped patterns is inspired from Gabor function. This function is used in Gabor filters, which is a linear filter used for texture analysis. This essentially means that it detects if there is any certain frequency content in an image along various directions in a confined region. The frequency and orientation representations of Gabor filters are compatible to those of a mammalian visual system \cite{marcelja80}. Because of its capacity to extract both spatial and frequency information from images, the Gabor filter is a prominent tool in computer vision and image processing.

Gabor filters, inspired by the physiological characteristics of the mammalian visual system, particularly the primary visual cortex (V1), emulate the receptive fields of neurons in V1. Neurons in V1 respond to specific spatial frequencies and orientations of visual stimuli. Gabor filters, as mathematical models, selectively respond to spatial frequencies and orientations, making them adept at capturing information about the texture's characteristics in images. They excel in tasks requiring discrimination between textures based on spatial frequency and orientation, finding applications in pattern recognition (e.g., face and fingerprint recognition), texture classification, and medical imaging for identifying abnormalities or specific structures \cite{sultan2020}. The appropriateness of Gabor filters for texture analysis lies in their ability to replicate the receptive fields of neurons in the mammalian visual cortex, making them valuable for tasks where capturing spatial frequency and orientation information is crucial, such as in medical image processing.


In Fig. \ref{pat.match.sine}, the pattern-matching technique is illustrated, where two groups are fitted together. To measure the stripe orientation with respect to the crossing angle bisector, we initially rotated the coordinates, aligning the bisector with the new $x$-axis. Consequently, the expected orientation $\gamma$ from pattern matching became 90\degree. This precisely is the bisector-normal hypothesis (see Sect. \ref{res.dis}), that was first hypothesised in \cite{naka1977} and then experimentally established in \cite{mullick_ploscb}.

\subsection{Pattern-Matching Technique 2.0}
\label{pat.match2.0}

In the earlier version of the pattern-matching technique, the pedestrian coordinates were fitted using a 2D sinusoid, by maximising $C$ given by Equation \ref{c.eq}, using the Nelder-Mead simplex algorithm. Out of several factors of merit, the objective function $C$ bears the importance of being a major benchmark tool, from an emergent patterns perspective, to gauge the realism of simulated data. When pedestrian coordinates from the two groups are fitted together, the maximum possible value of $C$ is $C_{max}=2$. A higher value of $C$, or a value of $C/C_{max}$ close to $1$, would signify better fitting, and consequently a better analysis of the results. However, for the previous computational method it was obtained that the median values of $C/C_{max}$ are in the range $0.45$ to $0.7$ \cite{mullick_ploscb}. So here we attempt to improve the fitting procedure by using different strategies.

\begin{figure}[h!]
    \centering
    \includegraphics[width=\textwidth]{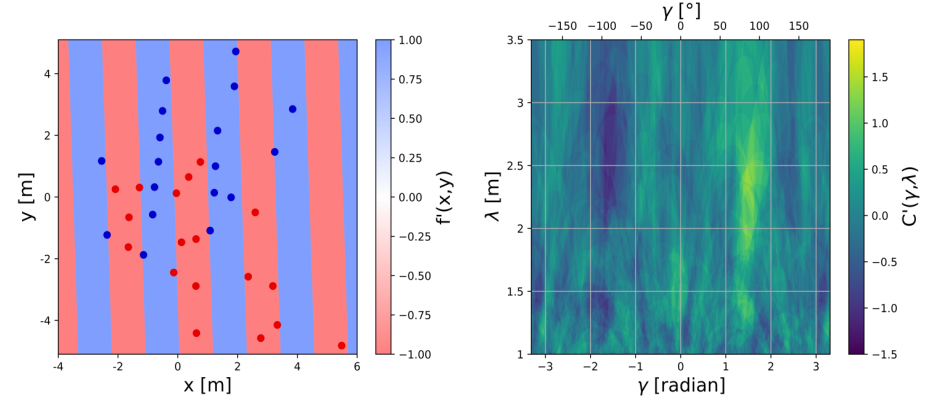}
    \caption{Demonstration of pattern-matching technique with the square wave for the typical trial shown in Figure \ref{pat.match.sine}. The left panel shows the 2D square wave $f'$ fitted on pedestrian coordinates from the two groups, denoted by blue and red dots, and the right panel shows the corresponding objective function $C'$ as a function of $\gamma$ and $\lambda$.}
    \label{pat.match.sign}
\end{figure}

In our modified pattern matching technique presented in this paper, we use a 2D square-wave, instead of a 2D sinusoid to fit the pedestrian coordinates. The 2D sine wave has varying values along its oscillations, representing the positions of pedestrians at different points where the stripe was the most visible. However, in our modified method, the sine wave is replaced by a square wave, which simplifies the pattern by assigning fixed values to the troughs and crests. Denoting the square wave by $f'$, its functional form could be given by \begin{equation}
    f'=\text{sgn}\Big[ \sin\Big(\frac{2\pi X}{\lambda}+\psi\Big)\Big],
    \label{gabor_sign}
\end{equation} where $X$ is already defined by (2). In the above expression, sgn is the signum function. Typically, the troughs are assigned a value of $1$, where we fit one of the groups of pedestrians, while the crests are assigned a value of $-1$, where the other group of pedestrians are fitted. The corresponding objective function, say $C'$, which is to be optimised now remains the same as (3) except for $f$ being replaced by $f'$; \begin{equation}
C' = \sum_{\text{group 1}}\frac{f'(x,y)}{N_1} + \sum_{\text{group 2}}\frac{-f'(x,y)}{N_2}.
\label{cp.eq}
\end{equation} Fig. \ref{pat.match.sign} shows the modified pattern-matching technique done using the square wave $f'$, and by maximising the objective function $C'$, done for the same trial shown in Fig. \ref{pat.match.sine}.

Another significant change that we introduce in the pattern matching technique is the optimisation algorithm. Here we use Simulated Annealing (SA) to maximise the objective function, instead of the previously used Nelder-Mead (NM) simplex algorithm. One of the vital shortcomings of NM algorithm is its tendency to converge to the local optima instead of the global optima. This behaviour arises because the algorithm relies on the given initial points of the parameters, and different initial points could lead to different final outcomes.

To address these challenges we apply the Simulated Annealing (SA) algorithm to the pattern-matching technique. Simulated Annealing is a stochastic optimization algorithm used to find approximate solutions to complex optimization problems. Inspired by the annealing process in metallurgy, it starts with an initial solution, and then iteratively explores the solution space by accepting both better and worse solutions with a probability determined by a temperature parameter. The temperature decreases over time, allowing the algorithm to make transition from a more explorative phase to a more exploitative one. This method is particularly effective for combinatorial optimization problems and situations where finding the global optimum is challenging due to the presence of multiple local optima.

By using simulated annealing, we aim to obtain more robust and reliable results for the optimal values of $\gamma$ and $\lambda$. SA explores the parameter space in a probabilistic manner, gradually reducing the search radius and focusing on promising regions. This approach can potentially overcome the limitations of the NM algorithm and provide better convergence to the global optima.


So to summarize our methods, we have used two different wave forms, viz. 2D sinusoid and 2D square wave, to capture the periodicity of the striped patterns, and for each of these wave forms we have used two different optimisation algorithms, viz. Nelder-Mead simplex and Simulated Annealing, to maximise the objective function for fitting pedestrian coordinates. Therefore, the combination of these strategies produce 4 different sets of pattern-matching results, which we are going to present and analyse in the next section.

\section{Results and Discussions}
\label{res.dis}

The primary objective of this paper is to analyse how pedestrian coordinates can be fitted, as accurately as possible, on the selected waveforms by using different optimisation algorithms. In this section we present the results obtained from the 4 different methods. To examine the fitting results and to compare the outcomes from various strategies defined in the previous section, we have also used statistical hypothesis tests. The comprehensive analysis provides further insights into the performance and reliability of the optimisation algorithms in the context of crowd motion studies.

\begin{figure}[h!]
    \centering
    \includegraphics[width=0.6\textwidth]{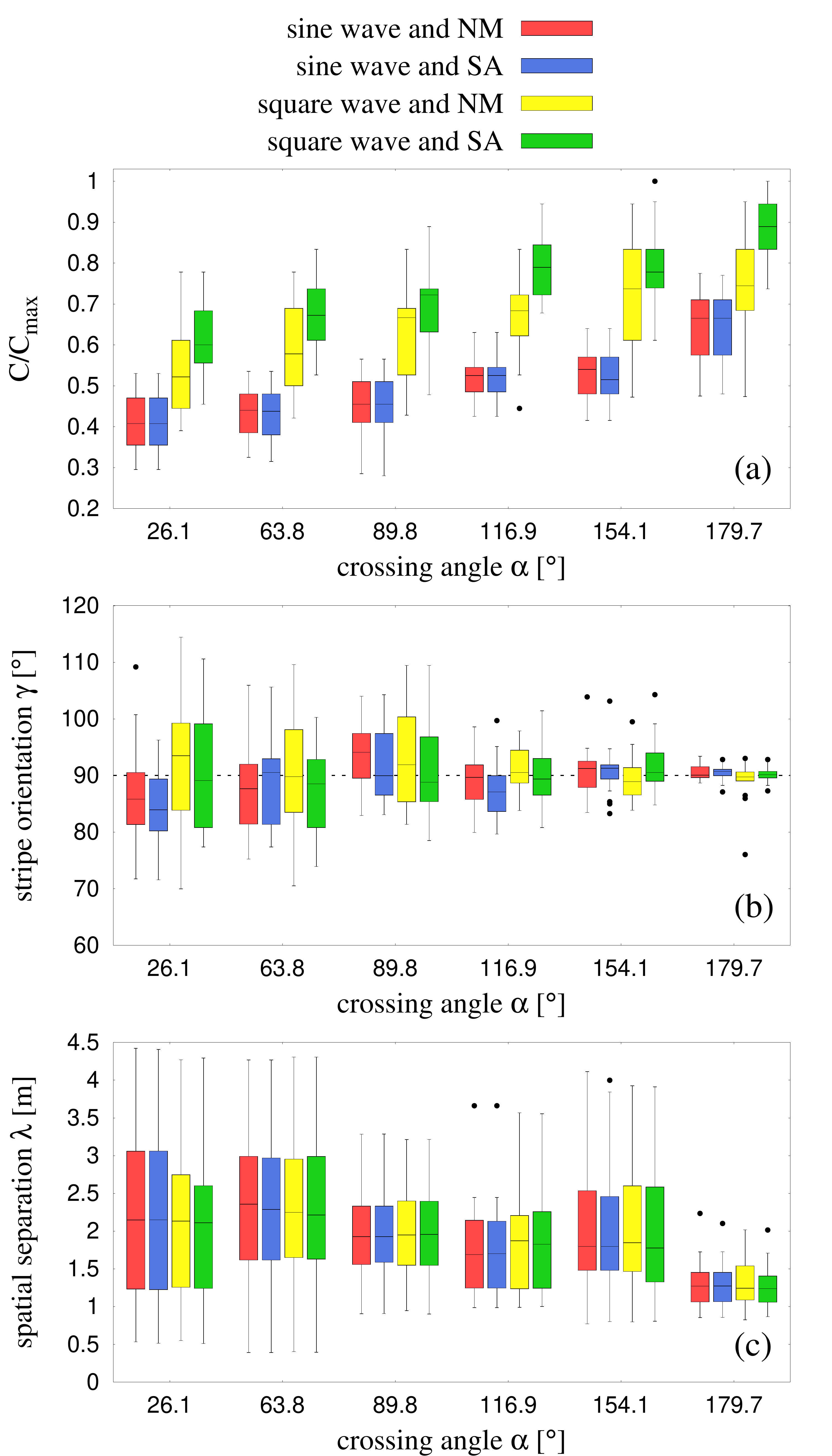}
    \caption{Boxplots to demonstrate the comparison of 4 sets of values of (a) $C/C_{max}$ (b) $\gamma$ and (c) $\lambda$ values for all the experimental trials, obtained by various strategies to implement the pattern-matching technique. In (b) the dashed black line indicates the $\gamma=90\degree$, which is the expected orientation of the striped patterns according to the bisector normal hypothesis.}
    \label{boxplots}
\end{figure}

In Figure \ref{boxplots} we present boxplots for the obtained values of the maximised (and normalised) objective function $C$ (or $C'$), stripe orientation $\gamma$ and spatial separation $\lambda$ estimated from 4 different optimisation strategies for 106 experimental trials. During the application of the pattern-matching technique using the Nelder-Mead (NM) simplex algorithm or Simulated Annealing (SA) with 2D sinusoid or 2D square wave, we obtained optimal values for $\lambda$ and $\gamma$ that maximize our objective function $C$ or $C'$. These optimal values from SA are very similar to the ones obtained using NM. This suggests that both algorithms are capable of finding the global optimum for our pattern-matching problem using either of the waveforms.

Let us first discuss the effect of changing the waveform from a 2D sinusoid to a 2D square wave. The transition from a continuously varying sine wave to a discontinuous square wave is clearly visualized in the plots for $C$ (Fig. \ref{pat.match.sine}) and $C'$ (Fig. \ref{pat.match.sign}), where we demonstrate pattern-matching (both using NM) for a typical experimental trial. The maximum value of $C'$ obtained for this trial is 1.888, whereas when we used the sine wave for pattern-fitting, the maximum value of $C$ was only 1.205. This indicates an improved fit of the pedestrians onto the troughs and crests of the square wave, effectively segregating the two groups. This is precisely the improvement that we intended to bring about by using a square wave instead of a sine wave for pattern matching. Boxplots for $C/C_{max}$ values for all the crossing angles in Fig. \ref{boxplots}(a) show that this improvement is systematic and could be seen for all the crossing angles of our data set.

It is also evident from Figure \ref{boxplots}(a) that for the sinusoidal wave, both NM and SA yield similar values of the maximised objective function. This has been verified by a series of one-way ANOVAs (see Table \ref{cvalues.stat.test}), where all $p$-values are $>0.8$ and values of $\eta^2$ are $\lesssim 10^{-3}$. However, for the square curve a significant increase in the maximised objective function is observed when switching to SA from NM. One-way ANOVAs (see Table \ref{cvalues.stat.test}) confirm this statistical significance (all $p$-values $<0.05$ and values of $\eta^2>0.1$). So undoubtedly, SA outperforms NM in achieving a more accurate fit of the pedestrians on the stripes, and this effect is more prominent while using a square curve compared to a sinusoid.

\begin{table}[h!]
    \centering
    \caption{Results of one-way ANOVAs to check statistical significance between maximised $C$ values obtained by different strategies of optimisation used in this research. The $F$-statistic in a one-way ANOVA assesses whether the means of more than two groups are significantly different; it is calculated by comparing the variance between group means to the variance within groups. The $p$-value associated with the $F$-statistic is the probability of observing an $F$-statistic as extreme as, or more extreme than, the one calculated from the sample data, assuming that there is no difference in population means. Asterisks indicate statistically significant $p-$values, for a chosen significance level of $0.05$. $\eta^2$ is a measure of effect size that indicates the proportion of variance that can be explained.}
    \begin{tabular}{|c|c|c|c|c|}
         \hline
         crossing angle&strategies tested&$F-$statistic&$p-$value&$\eta^2$\\
         \hline
         \multirow{2}{*}{$30\degree$}&sine+NM, sine+SA&$F(1,34)=0.004$&$0.953$&$\sim10^{-4}$\\
         \cline{2-5}
         &square+NM, square+SA&$F(1,34)=5.35$&$0.027*$&$0.136$\\
         \hline
         \multirow{2}{*}{$60\degree$}&sine+NM, sine+SA&$F(1,34)=0.065$&$0.801$&$0.002$\\
         \cline{2-5}
         &square+NM, square+SA&$F(1,34)=6.082$&$0.019*$&$0.152$\\
         \hline
         \multirow{2}{*}{$90\degree$}&sine+NM, sine+SA&$F(1,36)=0.001$&$0.982$&$\sim10^{-5}$\\
         \cline{2-5}
         &square+NM, square+SA&$F(1,36)=4.473$&$0.041*$&$0.11$\\
         \hline
         \multirow{2}{*}{$120\degree$}&sine+NM, sine+SA&$F(1,32)=0$&$0.988$&$<10^{-5}$\\
         \cline{2-5}
         &square+NM, square+SA&$F(1,32)=12.67$&$0.001*$&$0.284$\\
         \hline
         \multirow{2}{*}{$150\degree$}&sine+NM, sine+SA&$F(1,32)=0.024$&$0.878$&$\sim10^{-4}$\\
         \cline{2-5}
         &square+NM, square+SA&$F(1,32)=2.203$&$0.148$&$0.064$\\
         \hline
         \multirow{2}{*}{$180\degree$}&sine+NM, sine+SA&$F(1,32)=0$&$0.985$&$\sim10^{-5}$\\
         \cline{2-5}
         &square+NM, square+SA&$F(1,32)=11.88$&$0.002*$&$0.271$\\
         \hline
    \end{tabular}
    \label{cvalues.stat.test}
\end{table}

As mentioned in Sect. \ref{pat.match2.0}, NM as an optimisation algorithm has it’s disadvantages since it highly depends on the initial points for the parameters of the objective function. To overcome this limitation and obtain reliable results, we had to perform repetitive trials and optimisations using various initial points. The process of repetitive trials and optimisations significantly increased the computational time and required extensive manual intervention. It involved running the process multiple times with different initial parameter values and examining the results to select the best outcome. This labor-intensive procedure posed practical limitations, especially when dealing with large datasets and a complex optimisation problem, such as the pattern-matching.

\begin{figure}
    \centering
    \includegraphics[width=0.6\textwidth]{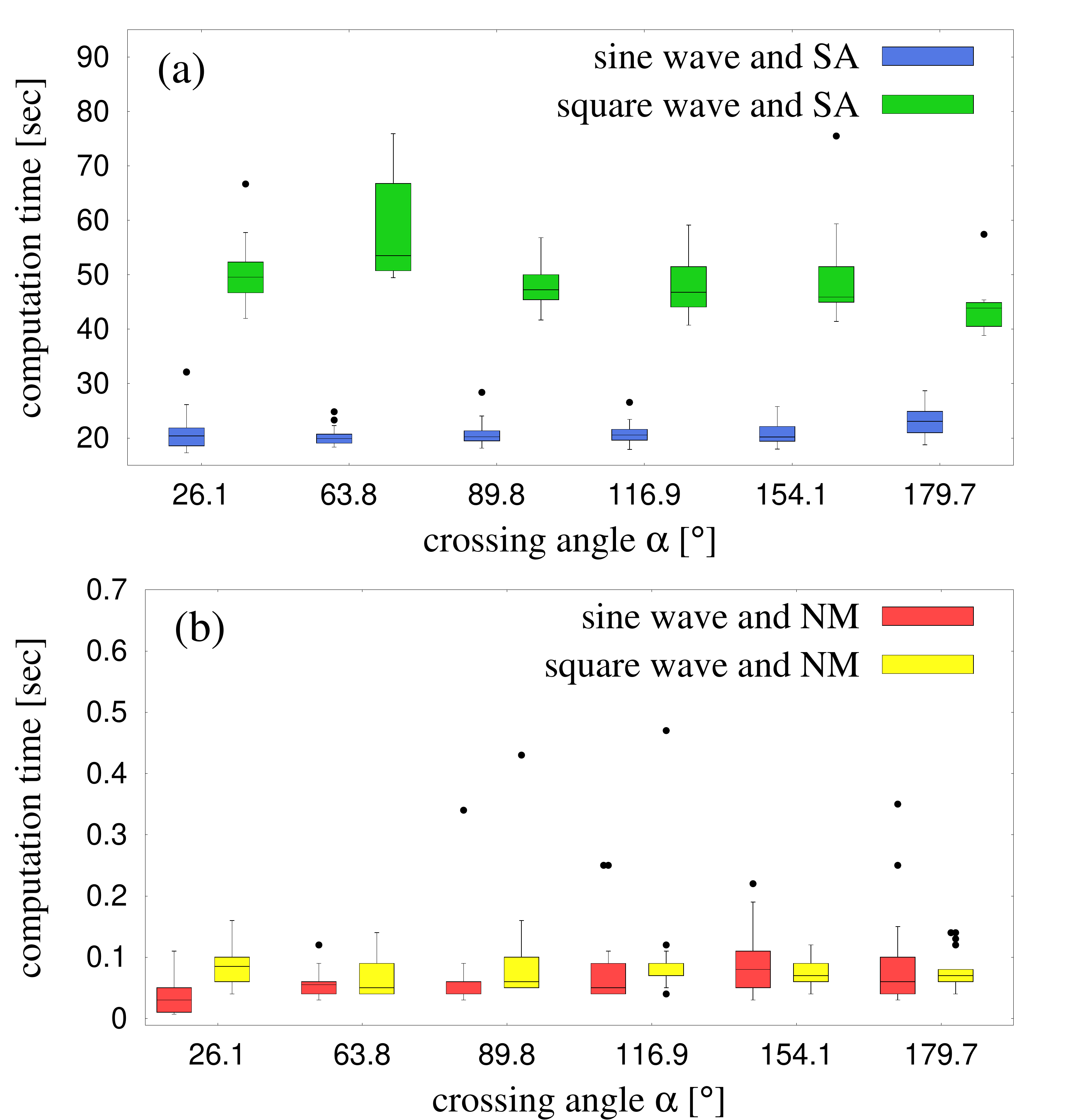}
    \caption{Boxplots to demonstrate the comparison of computation times required per trial of optimisation using (a) Simulated Annealing (SA) and (b) Nelder-Mead (NM) simplex algorithm. The computation time for SA is two order of magnitudes higher compared to that for NM.}
    \label{comp_time}
\end{figure}

One noticeable difference between the two algorithms was the computational time required. While pattern matching with SA took significantly longer time per single trial compared to that with NM, as shown by boxplots in Figure \ref{comp_time}, it was able to converge to the global maxima of the objective function in almost one iteration for all cases. In contrast, the NM required repetitive optimisation runs with different initial points for the parameters, which significantly increased the total computation time. The required numbers of repetitive runs depends solely on the quality of guesswork done to put the initial points. For a multi-parameter optimisation problem like ours, there is no wiser strategy to choose these initial points.

In Figures \ref{boxplots}(b) and \ref{boxplots}(c), we present the boxplots for the obtained values of stripe orientation $\gamma$ and inter-stripe spatial separation $\lambda$. For both the cases, the results do not appear to depend on the strategy of pattern-matching used. This, in a way, validates the accuracy of our extension of pattern-matching with different waveforms and optimisation algorithms.

A significant aspect of the stripe orientation $\gamma$ values is the validity of bisector-normal hypothesis. This hypothesis states that the self-organising stripes in crossing flows of pedestrians are perpendicular to the bisector of crossing angle, as depicted schematically in Figure \ref{stripe_schematic}. The validity of this hypothesis was established in \cite{mullick_ploscb}, where pattern-matching was performed using a 2D sinusoid and the Nelder-Mead simplex algorithm. Here we present 3 additional strategies of the pattern-matching procedure, and from Figure \ref{boxplots}(b) it seems that all the $\gamma$ values are distributed around $90\degree$. To verify this numerically we performed a series of one-sample $t$-tests to compare the means of each of these distributions with the hypothetical value $90\degree$. The results are summarised in Table \ref{t.tests}, which shows that when we used a square wave for pattern-matching, the obtained values of stripe orientation are not statistically different from $90\degree$. But when we use a sine wave for the pattern fitting we observe only one case each for NM and SA algorithms where the result is statistically significant. This perhaps once again establishes the fact, albeit only up to some extent in this case, that using a square wave instead of a sine wave improves the quality of results obtained from the pattern-matching procedure. 

\begin{table}[h!]
    \centering
    \caption{Results of one-sample $t$-tests performed for stripe orientation $\gamma$ values obtained from the 4 strategies of pattern-matching technique. The $t$-statistic in a one-sample $t$-test measures how far the sample mean deviates from the hypothesized population mean in terms of standard errors. On the other hand, the $p$-value is the probability of observing a $t$-statistic as extreme as, or more extreme than, the one calculated from the sample data, assuming that the true population mean is equal to the hypothesized mean. Asterisks indicate statistically significant $p$-values, where our chosen level of significance is $0.05$.}
    \begin{tabular}{|c|c|c|c|c|c|c|c|c|}
         \hline
         \multirow{2}{*}{crossing angle}&\multicolumn{2}{c|}{sine+NM}&\multicolumn{2}{c|}{sine+SA}&\multicolumn{2}{c|}{square+NM}&\multicolumn{2}{c|}{square+SA}\\
         \cline{2-9}
         &$t$-statistic&$p$-value&$t$-statistic&$p$-value&$t$-statistic&$p$-value&$t$-statistic&$p$-value\\
         \hline
         $30\degree$&$t(17)=-1.284$&$0.216$&$t(17)=-3.625$&0.002*&$t(17)=1.14$&0.27&$t(17)=-0.035$&0.972\\
         \hline
         $60\degree$&$t(17)=-0.825$&0.421&$t(17)=-0.234$&0.818&$t(17)=0.218$&0.83&$t(17)=-1.55$&0.139\\
         \hline
         $90\degree$&$t(18)=3.092$&0.006*&$t(18)=1.511$&0.148&$t(18)=1.581$&0.131&$t(18)=0.597$&0.558\\
         \hline
         $120\degree$&$t(16)=-0.804$&0.433&$t(16)=-2.204$&0.042&$t(16)=0.627$&0.542&$t(16)=0.056$&0.956\\
         \hline
         $150\degree$&$t(15)=0.534$&0.6&$t(16)=0.782$&0.445&$t(16)=-0.76$&0.458&$t(16)=1.49$&0.156\\
         \hline
         $180\degree$&$t(16)=1.747$&0.01&$t(16)=1.693$&0.11&$t(16)=-0.941$&0.361&$t(16)=0.639$&0.532\\
         \hline
    \end{tabular}
    \label{t.tests}
\end{table}


\section{Conclusions}
\label{con}

We considered an experimental data for crossing flows of two groups of pedestrians, whose trajectories were recorded while they attempted to cross each other. The numerical strategies described in this paper were utilised to detect the self-organising striped patterns in the crossing region of the two flows and to study their geometric properties. Previous research employed a method known as the pattern-matching technique which basically is fitting a 2D sinusoid to pedestrian coordinates using the Nelder-Mead (NM) simplex algorithm. We extended the pattern-matching method by using a 2D square wave to fit the data, and switching to Simulated Annealing (SA) as the optimisation algorithm.

In this research our aim was to get an improved fit of the pedestrians on the striped pattern. We expected that switching from a sine wave to a square wave would definitely bring about the desired improvement. However, SA in combination with the square wave demonstrated a further enhancement in the quality of pattern-fitting, which was obvious by no means. The advantage of SA was its ability to explore the parameter space more effectively and consistently in order to reach the global optima. On the other hand, NM's dependence on the initial points made it more prone to converging to local optima and necessitated additional manual intervention to ensure reliable results. Overall, SA performed better than NM for both the waveforms used in this research, as shown in Figure \ref{boxplots}(a) and confirmed by one-way ANOVAs.

Gaining a deep understanding of the behaviour exhibited by large human crowds is essential for effective crowd management planning. Computational modeling and simulation of crowds play a pivotal role in this endeavour. The objective function $C$ (or $C'$), which is at the core of pattern-matching technique, has the potential to act as a benchmark to judge the quality of artificially simulated trajectories, as in how realistic they are. The key concept would be to verify the presence and geometrical properties of the self-organising patterns. Improvements in pattern-matching presented in this paper would aid in the understanding of complex dynamics of crowds through numerical simulation and eventually to develop effective crowd management strategies.

In conclusion, while our primary focus has been on the dynamics of pedestrian flows, the versatility of our pattern matching method extends its applicability to diverse contexts. Beyond pedestrian dynamics, the underlying principles of our method render it adaptable to numerous practical scenarios where the identification of self-organizing patterns holds scientific significance. For example: (1) analyzing the flow of vehicles in transportation systems, such as highways, intersections, or parking lots, to optimize traffic management and improve overall efficiency; (2) understanding patterns in the movement of cyclists, both in urban environments and dedicated cycling lanes, to enhance cycling infrastructure and safety; (3) observing and analyzing the movement patterns of animals in their natural habitats for ecological studies and wildlife conservation; (4) studying the movement patterns of robots or automated systems in manufacturing plants to enhance efficiency, minimize collisions, and improve overall automation processes; (5) applying the methodology to study the movement patterns of ocean currents, which can have implications for marine navigation, environmental monitoring, and climate studies; (6) investigating the movement patterns of aircraft in airport taxiways or on runways to optimize air traffic control and ground operations; (7) observing the movement of components or products along assembly lines in manufacturing facilities to enhance production efficiency and identify potential bottlenecks. These examples highlight the versatility of our method, showcasing its potential in various domains where understanding and optimizing movement patterns are essential.\\[1cm]





%
%

\noindent \textbf{Declarations}\\

\small\noindent \textbf{Acknowledgements}

\noindent Fruitful discussions with C. Appert-Rolland, J. Pettré and W. H. Warren are greatly acknowledged, which lead to the inception of this problem.\\

\small\noindent \textbf{Author contributions}

\noindent Conceptualization: PM, methodology: SW and PM, software: SW and PM, validation: SW and PM, formal analysis: SW and PM, investigation: SW, resources: PM, data curation: SW, writing -- original draft preparation: SW, writing -- review and editing: PM, visualization: SW and PM, supervision: PM, project administration: PM, funding acquisition: PM. All authors have read and agreed to the current version of the manuscript.\\

\small\noindent \textbf{Funding}

\noindent PM acknowledges financial support from National Science Center (NCN, Poland) through SONATA grant no. 2022/47/D/HS4/02576.\\

\small\noindent \textbf{Availability of data and materials}

\noindent The computational codes (in R) used for pattern-matching technique and the outputs are available at: \url{https://doi.org/10.5281/zenodo.8348304}. The trajectories dataset for crossing flows are available at: \url{https://doi.org/10.5281/zenodo.5718430}.\\

\small\noindent \textbf{Competing interests}

\noindent The authors declare that there are no competing interests.\\

\bibliography{bibliography}   


\end{document}